\documentclass[a4paper,12pt]{amsart}
\usepackage{ifthen}
\usepackage{graphicx}
\usepackage{mathrsfs}
\usepackage{color}
\nonstopmode
\numberwithin{equation}{section}
\newtheorem{thm}{Theorem}[section]

\newtheorem{lem}[thm]{Lemma}
\newtheorem{prop}[thm]{Proposition}

\theoremstyle{definition}

\newtheorem{example}[thm]{Example}

\newenvironment{pf}[1][]{%
 \vskip 3mm
 \noindent
 \ifthenelse{\equal{#1}{}}%
  {{\slshape Proof. }}%
  {{\slshape #1.} }%
 }%
{\qed\bigskip}

\newcommand{\C}{{\mathbb C}}
\newcommand{\D}{{\mathbb D}}

\newcommand{\rhp}{{H}}

\newcommand{\sh}{{\mathcal S}_\mathrm{H}}

\newcommand{\hol}{{\operatorname{Hol}}}
\newcommand{\HT}{{\mathcal{HT}}}
\newcommand{\har}{{\operatorname{Har}}}

\newcommand{\Li}{{\,\operatorname{Li}}}

\newcommand{\cm}{{\mathcal T}}
\newcommand{\acm}{{\widetilde{\mathcal T}}}

\renewcommand{\Im}{{\,\operatorname{Im}\,}}

\renewcommand{\Re}{{\,\operatorname{Re}\,}}

\newcommand{\Gauss}{{\null_2F_1}}

\newcommand{\inv}{^{-1}}

\newcommand{\term}[1]{{\it #1}}

\newcommand{\aand}{{\quad\text{and}\quad}}

\newcounter{minutes}
\divide\time by 60
\newcounter{hours}
\multiply\time by 60 \addtocounter{minutes}{-\time}

\begin{document}
\bibliographystyle{amsplain}
\title
{Completely monotone sequences and harmonic mappings }
\author[B.-Y. Long]{Bo-Yong Long}
\address{School of Mathematical Sciences, Anhui University, Hefei  230601, China}
\email{boyonglong@163.com}
\author[T. Sugawa]{Toshiyuki Sugawa}
\address{Graduate School of Information Sciences \\
Tohoku University\\
Aoba-ku, Sendai 980-8579, Japan} \email{sugawa@math.is.tohoku.ac.jp}
\author[Q.-H. Wang]{Qi-Han Wang}
\address{School of Mathematical Sciences, Anhui University, Hefei  230601, China}
\email{qihan@ahu.edu.cn}
\begin{abstract}
In the present paper, we will study geometric properties of harmonic mappings
whose analytic and co-analytic parts are (shifted) generated functions of completely monotone  sequences.
\end{abstract}

\thanks{
The present research was supported in part by
Natural Science Foundation of Anhui Province (1908085MA18),
Foundation of Anhui Educational Committee (KJ2020A0002), China.
}

\maketitle

\section{Introduction and preliminaries}

A sequence $\{a_{n}\}_{n=0}^\infty$ of real numbers  is called  \term{completely monotone}
(or \term{totally monotone}) if
$\Delta^k a_n\ge 0$ for all integers $n,k\ge0.$
Here, $\Delta^ka_n$ is defined recursively by $\Delta^0 a_n=a_n, n\ge 0,$ and
$$
\Delta^{k} a_{n}:=\Delta^{k-1} a_{n}-\Delta^{k-1} a_{n+1},  \quad n\geq 0, \quad k\geq 1.
$$
Note that a completely monotone sequence is non-negative, non-increasing and convex.
Hausdorff  \cite{Haus21} showed that  $\{a_{n}\}$  is a completely monotone
sequence precisely when there is a positive Borel measure $\mu$ on $[0, 1]$ such that
\begin{align}\label{Hausdorff1}
a_{n}=\int^{1}_{0}t^{n}d\mu(t), \quad n\geq 0.
\end{align}
Therefore, the word ``completely monotone sequence" is a synonym of ``Hausdorff moment sequence".
When $a_0=1$ the sequence is said to be \term{normalized}.
Note that the condition $a_0=1$ means that $\mu$ is a probability measure.
In particular, the generating function of the normalized Hausdorff sequence $\{a_n\}$ is represented in the form
\begin{align}\label{eq:Haus2}
F(z)=1+\sum_{n=1}^{\infty}a_{n}z^{n}=\int^{1}_{0}\frac{d\mu(t)}{1-tz}
\end{align}
for a Borel probability measure $\mu$ on $[0,1].$
We denote by $\cm$ the set of those functions $F$ generated by normalized Hausdorff moment sequences.
For instance, letting $\mu$ be the Dirac measure with unit mass at $t=0$ or $t=1,$ we see that
the functions $F_0(z)=1$ and $F_1(z)=1/(1-z)$ belong to $\cm.$
By the form \eqref{eq:Haus2}, we observe that a function $F\in\cm$
is analytically continued to the slit domain $\Lambda:=\mathbb{C}\backslash[1,+\infty)$.
We note that $F(x)$ is non-decreasing in $-\infty<x<1$ and, in particular, that $F(x)$ has a limit
(possibly $+\infty$) as $x\to 1^-.$
The value of the limit will be denoted by $F(1^-).$
By the form of $F,$ we observe that
\begin{equation}\label{eq:F}
|F(z)|\le \sum_{n=0}^\infty a_n|z|^n=F(|z|)\le F(1^-)=\sum_{n=0}^\infty a_n,
\quad |z|<1.
\end{equation}
Note also that $F(1^-)\ge a_0=1.$
We denote by $\acm$ the set of \term{shifted} generated functions $zF(z)$ for $F\in\cm.$
For instance, the functions $f_0(z)=z$ and $f_1(z)=z/(1-z)$ both are members of $\acm.$

Completely monotone sequences are closely related with moment problems and the theory of continued fractions
and thus important not only in analysis but also in probability and applied mathematics,
see \cite{LP16, ST:moment, Wall:fraction} for instance.
In recent years, the theory of universally prestarlike functions (containing universally convex and universally starlike functions)
was developed and  an intimate connection with $\cm$ was found (see \cite{Rusch07} and \cite{RSS09}).
As is well recognized, many kinds of special functions may be described in terms of functions in $\cm$
(see \cite{BRS18}, \cite{RSS09} as well as Section 4 below).

Let $\hol(\Lambda)$ denote the set of analytic functions on the domain $\Lambda=\C\setminus[1,+\infty).$
The following lemma is more or less known to experts.
This sort of result was formulated in \cite[Lemma 2.1]{RSS09} and then simplified
by Liu and Pego \cite{LP16} (see Remark 2 therein).

\begin{lem}\label{lem:char}
Let $F\in \hol(\Lambda)$.
Then $F\in\cm$, i.e.  $F$ can be represented in the form
$$F(z)=\int^{1}_{0}\frac{d\mu(t)}{1-tz} $$
for a Borel probability measure $\mu$ on $[0,1]$, if and only if the following three conditions are fulfilled:
\begin{enumerate}
\item[(i)] $F(0)=1$;
\item[(ii)] $F(x)$ is a non-negative real number for each $x\in(-\infty,1)$;
\item[(iii)] $\Im F(z)\geq 0$ whenever $\Im z>0$.
\end{enumerate}
Moreover, the measures $\mu$ and the functions $F$ are in one-to-one correspondence.
\end{lem}

Let $\har(\D)$ denote the class of complex-valued harmonic functions on the unit disk $\D.$
Then, each function $f$ in $\har(\D)$ is uniquely expanded in the form
\begin{equation}\label{eq:harm}
f(z)=\sum_{n=-\infty}^\infty a_n r^{|n|}e^{i n\theta},  \quad z=re^{i\theta}\in\D.
\end{equation}
For another $F(z)=\sum_n A_n r^{|n|}e^{i n\theta}$ in $\har(\D),$
we define the (harmonic) \term{convolution} (or the \term{Hadamard product})
of $f$ and $F$ by
$$
(f\ast F)(z)=\sum_{n=-\infty}^\infty a_n A_n r^{|n|}e^{i n\theta}.
$$
Note that $f\ast F\in\har(\D)$ whenever $f,F\in\har(\D).$
It is often more convenient to express $f$ in \eqref{eq:harm} in the form
$$
f(z)=\sum_{n=1}^\infty a_{-n}\bar z^n+\sum_{n=0}^\infty a_nz^n
=\overline{g(z)}+h(z),
$$
where
$$
h(z)=\sum_{n=0}^\infty a_nz^n
\aand
g(z)=\sum_{n=1}^\infty b_nz^n
\quad (b_n=\overline{a_{-n}},~ n\ge 1).
$$
The analytic functions $h$ and $g$ are called the \term{analytic part} and
the \term{co-analytic part} of $f,$ respectively.
The convolution of harmonic functions $f_1=h_1+\bar g_1$ and $f_2=h_2+\bar g_2$ is described also by
$$
f_1\ast f_2=h_1\ast h_2+\overline{g_1\ast g_2},
$$
where $h_1\ast h_2$ and $g_1\ast g_2$ are the ordinary Hadamard products.
See \cite{Rus:conv} for basics of convolutions of anaytic functions.

A smooth map $f:\D\to\C$ is locally univalent at $z_0$ if
the Jacobian $J_f=|f_z|^2-|f_{\bar z}|^2$ does not vanish at $z_0$
by the Inverse Mapping Theorem.
Lewy's theorem asserts that the converse is true for harmonic mappings.
Therefore, a harmonic mapping $f=h+\bar g$ is locally univalent and
sense-preserving at $z$
if and only if $J_f(z)=|f_z(z)|^2-|f_{\bar z}(z)|^2=|h'(z)|^2-|g'(z)|^2>0.$
In particular, then we have $h'(z)\ne 0$ and the function
$$
\omega_f(z)=\frac{g'(z)}{h'(z)}
$$
is holomorphic at $z$ and satisfies the inequality $|\omega_f(z)|<1.$

We denote by $\sh$ the set of sense-preserving harmonic \term{univalent} functions $f$ in
$\har(\D)$ normalized by $f(0)=f_z(0)-1=0.$
In what follows, we will mean sense-preserving and injective (one-to-one)
by the term ``univalent".
Set also $\sh^0=\{f\in\sh: f_{\bar z}(0)=0\}.$
These classes were introduced and studied by Clunie and Sheil-Small \cite{CS84}.
Nowadays, many researchers are studying them and their subclasses intensively.
See the monograph \cite{Duren:harm} for fundamental theory and recent progress
of harmonic univalent mappings.
Note that $\omega_f=\bar f_z/f_z=g'/h'$ satisfies the inequality $|\omega_f|<1$
on $\D$ for $f=h+\bar g\in\sh.$
The quantity $\omega_f=\bar f_z/f_z=\overline{f_{\bar z}}/f_z$ is called
the \term{second complex dilatation} of $f.$
If $f$ is univalent and if $|\omega_f|\le k$ for a constant $k<1,$
the mapping $f$ is called $k$-quasiconformal (or $K$-quasiconformal
in the most of the literature, where $K=(1+k)/(1-k)$).
For the theory of quasiconformal mappings, the reader should consult
the standard monograph \cite{Ahlfors:qc2} by Ahlfors.
Much attention has been paid to the class of harmonic quasiconformal mappings
on the unit disk.
See \cite{CP19, Kalaj10, PSZ18, SRJ16} and references therein.

For a constant $c$ with $|c|<1,$ we define the class
$\HT(c)$ to be the set of functions $f\in\har(\Lambda)$ of the form
$$
f(z)=h(z)+c\overline{g(z)}
$$
for some $h,g\in\acm.$
In other words, each member $f$ of $\HT(c)$ is represented as
$$
f(z)=\int_0^1\frac{z}{1-tz}d\mu(t)+c \int_0^1\frac{\overline{z}}{1-t\bar z}d\nu(t)
$$
for Borel probability measures $\mu,\nu$ on $[0,1].$
The purpose of this article is to study geometric properties of functions in $\HT(c)$
such as univalence, convexity in one direction, and quasiconformality.

In the next section, main results of this paper will be presented.
Their proofs are given in Section 3.
In Section 4, we will give a couple of examples and apply some of the main results
to polylogarithms and shifted hypergeometric functions.

\section{Main results}

\begin{thm}\label{Theorem 1.2}
Let $c$ be a real constant with $0\le c<1.$
For $f\in\HT(c)$ and a constant $a\ge0,$ the following inequality holds:
\begin{align*}
|a+f(z)|\geq a+f(-|z|)\ge a+\lim_{r\to1^-}f(-r),
\quad z\in\D.
\end{align*}
\end{thm}

We remark that the above inequality is meaningful only when $a+f(-|z|)>0.$
In particular, if $a+\lim_{r\to 1^-}f(-r)\ge 0$ and if $f$ is non-constant,
then we conclude that $a+f(z)$ is non-vanishing on $|z|<1.$
If $c=0$, this theorem reduces to the main result of \cite{Pey65}.

Let $\rhp=\{z: \Re z<1\}$. Then we get the following result.

\begin{thm} \label{Theorem 1.3}  Let $f=h+c\bar g\in\HT(c)$ for a real constant $c$ with $0\le c<1.$
Suppose $h(z)=\int_0^1z(1-tz)\inv d\mu(t)$ and $g(z)=\int_0^1z(1-tz)\inv d\nu(t)$
for Borel probability measures $\mu$ and $\nu$ on $[0,1].$
Then for $z=x+iy\in\rhp$ with $y\ne0,$ unless $f$ is a constant function, the following hold:
\begin{enumerate}
\item[(i)]
$y\dfrac{\partial}{\partial y}\Re f(z)<0$;
\medskip
\item[(ii)]
$y\dfrac{\partial}{\partial x}\Im f(z)>0$ provided that $\mu=c\nu+(1-c)\lambda$
for a Borel probability measure $\lambda$ on $[0,1].$
\end{enumerate}
\end{thm}

Wirths \cite{Wirths75} proved the following useful result.

\begin{lem}\label{lem:Wirths}
Each function $h\in\acm$ is univalent on $\rhp=\{z: \Re z<1\}$
and the image domain $D=h(H)$ is convex in the direction of the imaginary axis.
\end{lem}

Here and hereafter, a domain $D$ in $\C$ is said to be
\term{convex in the direction of the imaginary axis}
if the intersection of $D$ with each line parallel to the imaginary axis is connected (or empty).
We can extend this result to the harmonic case.

\begin{thm} \label{Theorem 1.4}
Let $c$ be a real constant with $0\le c<1$ and $f\in\HT(c).$
If $f$ is locally univalent on $\rhp=\{z: \Re z<1\},$ then $f$ is univalent on $\rhp$ and
the image $f(\rhp)$ is convex in the direction of the imaginary axis.
\end{thm}

Unfortunately, we cannot drop the local univalence of $f$ in the assumption.
See Example \ref{ex:qc} in Section 4.

A linear combination is an important method to construct a new function, {\it cf.} \cite{LD18, SRJ16}.
However, it is well known that the convex combination of two univalent analytic functions
is not necessarily univalent, let alone convex combination of two univalent harmonic mappings.
The harmonic convolution $f_1\ast f_2$ of two harmonic functions $f_1$ and $f_2$ in $\har(\D)$
does not necessarily enjoy properties of $f_1$ or $f_2$, such as convexity or even (local) univalence
(see \cite{DNW12} for instance).
However, the following proposition shows that
the harmonic convolution and convex combinations keep the family $\HT(c)$ invariant in some sense.

\begin{prop}\label{proposition1.5}
Let $c_1, c_2$ be complex constants with $|c_j|<1.$
Then, for $f_j\in \HT(c_j),~j=1,2$, the following hold:
\begin{enumerate}
\item[(i)]
$sf_1+(1-s)f_2\in \HT(c)$ for $0\leq s\leq 1$ if $c=c_1=c_2;$
\item[(ii)]
$f_1\ast f_2\in \HT(c_1c_2).$
\end{enumerate}
\end{prop}


The next result gives us a sufficient condition for a function in $\HT(c)$ with a constant
$|c|<1$ to be quasiconformal on $\D.$

\begin{thm}\label{Theorem 1.6}
Let $h\in \acm$ and let $c$ be a real constant with $0\le c<1.$
Suppose that
\begin{align}\label{Theorem1.6.1}
\left|\frac{h'(tz)}{h'(z)}\right|\leq M,\quad z\in\D,~ 0\le t\le 1.
\end{align}
Then for each $g\in\acm$, the harmonic mapping
\begin{align}\label{equ1.12}
f(z):=h(z)+c\,\overline{(h\ast g)(z)},\quad z\in \mathbb{D}
\end{align}belongs to  $\HT(c)$ and is $k$-quasiconformal if $cM\le k<1$.
\end{thm}

Note that $h'(z)$ has no zeros on $\D$ since $h\in\acm$ is univalent on $\rhp$
by Lemma \ref{lem:Wirths}.
Letting $t=0,$ we observe that the condition $|h'(z)|\ge 1/M,~ z\in\D,$ is necessary
for \eqref{Theorem1.6.1}.
However, it is not easy to check \eqref{Theorem1.6.1} in general.
The following result may be helpful to find a value of $M.$

\begin{prop}\label{prop:M}
Let $h$ be an analytic function on $\D$ with $h'(0)=1$ and let $m$ be a positive constant.
Suppose that the following inequality holds:
\begin{equation}\label{eq:m}
\Re\frac{zh''(z)}{h'(z)}>-m,\quad z\in\D.
\end{equation}
Then
$$
\left|\frac{h'(tz)}{h'(z)}\right|\leq e^{2m},\quad z\in\D,~ 0\le t\le 1.
$$
\end{prop}

It is well known that a normalized analytic function $h(z)=z+a_2z^2+\dots$
maps $\D$ univalently onto a convex domain if and only if $\Re[zh''(z)/h'(z)]\ge -1.$
Therefore, we can take $1$ as the constant $m$ in the above proposition for this $h.$

Let $h, g\in \acm$.
It is important to look at
the quotient of the derivatives of two functions in $\acm$
when considering the local univalence or the quasiconformality of
the function of $\HT(c)$ for a consant $c$.
Under what conditions $g'/h'$ belong to $\cm$?
This question is interesting in itself.
In this context the following result proves to be useful.

\begin{thm}\label{Theorem 1.8}
Let $h, g\in \acm$ be represented by
\begin{align}\label{eq:hg}
h(z)=\int^{1}_{0}\frac{z\phi(t)}{1-tz}dt, \quad g(z)=\int^{1}_{0}\frac{z\psi(t)}{1-tz}dt
\end{align}
for nonnegative Borel functions $\phi$ and $\psi$ on $(0,1)$ with $\int_0^1\phi(t)dt=\int_0^1\psi(t)dt=1.$
If the inequality
\begin{align}\label{Theorem1.8.2}
\phi(s)\psi(t)\geq\phi(t)\psi(s)
\end{align}
holds for $0<s\leq t<1$, then $g/h$ and $g'/h'$ both belong to $\cm$.
\end{thm}

Note that the claim $g/h\in\cm$ was first proved in \cite[Theorem 1.10]{RSS09}.
When $\phi$ is non-vanishing, the condition in \eqref{Theorem1.8.2}
means that the function $\psi(t)/\phi(t)$ is non-decreasing in $0<t<1.$

Using Theorem \ref{Theorem 1.8}, we obtain another sufficient condition
for a function in $\HT(c)$ to be quasiconformal.

\begin{thm}\label{Theorem 1.9}
Under the hypotheses of Theorem \ref{Theorem 1.8}, further assume that
the function $F(z)=g'(z)/h'(z)$ has a finite limit $F(1^-).$
Then the function $f=h+c\bar g$ is $k$-quasiconformal on $\D$
if $0\le cF(1^-)\le k<1.$
\end{thm}

If $\lim_{x\to1^-}g'(x)=g'(1^-)<+\infty,$ then we have $F(1^-)=g'(1^-)/h'(1^-).$
If $g'(1^-)=+\infty,$ then $h'(1^-)=+\infty$ by the assumption $F(1^-)<+\infty.$
In this case, we may use l'H\^ospital's rule if the right-most limit below exists:
$$
F(1^-)
=\lim_{x\to1^-}\frac{g'(x)}{h'(x)}
=\lim_{x\to1^-}\frac{g''(x)}{h''(x)}.
$$

\section{Proofs of the main results}

In this section, we prove all the results in the previous section.

\begin{pf}[Proof of Theorem \ref{Theorem 1.2}]
By assumption, $f=h+c\overline{g}$ for some $g,h\in\acm.$
We first note the inequalities for $z=x+iy$ with $r=|z|<1$ and $0\le t\le 1:$
$$%
\Re\frac{z}{1-tz}=\frac{x-tr^2}{1-2tx+t^2r^2}\ge\frac{-r-tr^2}{1+2tr+t^2r^2}
=\frac{-r}{1+tr}\ge \frac{-1}{1+t}
$$%
because the function $x\mapsto (x-tr^2)/(1-2tx+t^2r^2)$ is increasing in $-r\le x\le r.$
Letting $\mu$ and $\nu$ be the representing measures of $h$ and $g,$ respectively,
we therefore have the estimates for $z$ with $|z|=r<1:$
\begin{align*}
|a+f(z)|&\ge a+\Re f(z)
=a+\Re h(z)+c \Re \overline{g(z)}
=a+\Re h(z)+c \Re g(z) \\
&=a+\int^{1}_{0}\Re\left(\frac{z}{1-tz}\right)d\mu(t)
+k\int^{1}_{0}\Re\left(\frac{z}{1-tz}\right)d\nu(t) \\
&\ge a+\int^{1}_{0}\frac{-r}{1+tr}d\mu(t)
+c\int^{1}_{0}\frac{-r}{1+tr}d\nu(t) \\
&=a+h(-r)+cg(-r)=a+h(-r)+c\overline{g(-r)}=a+f(-r) \\
&\ge a+\lim_{r\to1^-}f(-r).
\end{align*}
\end{pf}

\begin{pf}[Proof of Theorem \ref{Theorem 1.3}]
For $z=x+yi,$ by a straightforward computation, we have the expression
\begin{align*}
\Re f(z)&= \int^{1}_{0}\frac{x-t(x^2+y^2)}{1-2xt+t^{2}(x^2+y^2)}(d\mu(t)+cd\nu(t))
\aand \\
\Im f(z)&= \int^{1}_{0}\frac{y}{1-2xt+t^{2}(x^2+y^2)}(d\mu(t)-cd\nu(t)).
\end{align*}
Therefore, we have
\begin{align*}
\frac{\partial}{\partial y}\Re f(z)
&= \int^{1}_{0}\frac{-2yt(1-xt)}{(1-2xt+t^{2}(x^2+y^2))^2}(d\mu(t)+cd\nu(t)),
\aand\\   
 \frac{\partial}{\partial x}\Im f(z)
&= \int^{1}_{0}\frac{2yt(1-xt)}{(1-2xt+t^{2}(x^2+y^2))^2}(d\mu(t)-cd\nu(t)).
\end{align*}
Since $\mu+c\nu$ and $\mu-c\nu=(1-c)\lambda$ are positive measures,
we have the required inequalities for $x<1$ and $y\ne0.$
\end{pf}

For the proof of Theorem \ref{Theorem 1.4}, we need to recall the shear construction
developed by Clunie and Sheil-Small \cite{CS84}.
The following form is a vertical version of
a theorem of Clunie and Sheil-Small \cite[Theorem 5.3]{CS84}
(see also \cite[p.~37]{Duren:harm}) .

\begin{lem}[Clunie and Sheil-Small]\label{CS}
Let $f=h+\bar g$ be a locally univalent harmonic mapping on $\D.$
Then $f$ maps $\D$ univalently onto a convex domain in the direction of the imaginary axis
if and only if the analytic function $F=h+g$ maps $\D$ univalently onto
a convex domain in the direction of the imaginary axis.
\end{lem}

We denote by $f_*$ the $\pi/2$-rotation $if$ of $f$ about the origin.
Then $f=h+\bar g$ is convex in the direction of the imaginary axis
if and only if $f_*=if=ih+i\bar g=h_*-\overline {g_*}$ is convex in the direction
of the real axis.
Therefore, the above version follows from the original version \cite[Theorem 5.3]{CS84}.

\begin{pf}[Proof of Theorem \ref{Theorem 1.4}]
We will combine the technique employed by Wirths \cite{Wirths75} with the shear construction.
First note that the M\"obius transformation
$$
\varphi(\zeta)=\frac{2\zeta}{1+\zeta}
$$
maps $\D$ onto $H$.
Suppose that $f=h+c\bar g\in\HT(c)$ for some $0\le c<1$ and let
$F=h+cg.$
Then $F/(1+c)\in\acm$ and the proof of Wirth's theorem (Lemma \ref{lem:Wirths})
in \cite{Wirths75}
now implies that $F_1=F\circ\varphi=h\circ\varphi+cg\circ\varphi$ is univalent on $\D$
and convex in the direction of the imaginary axis.
Since $f_1:=f\circ\varphi=h\circ\varphi+\overline{cg\circ\varphi}$ is
locally univalent by assumption, now Lemma \ref{CS} implies that $f_1$ is univalent on $\D$
and convex in the direction of the imaginary axis.
Since $f=f_1\circ\varphi\inv,$ the assertion now follows.
\end{pf}

To prove the second part of Proposition \ref{proposition1.5}, we need the following lemma.

\begin{lem}\label{lem:conv}
Let $f, g\in \acm.$ Then $f\ast g\in\acm.$
\end{lem}

This fact is known to experts (see Roth, Ruscheweyh and Salinas \cite[p.~3172]{RRS08}).
Let us, however, give a proof because the authors could not find a proof in the literature.

\begin{pf}
Let $f(z)=\sum a_nz^{n+1}$ and $g(z)=\sum b_nz^{n+1}$
for normalized completely monotone sequences $\{a_n\}$ and $\{b_n\}.$
We have to show that $\{a_nb_n\}$ is completely monotone, too.
We first note the formula
\begin{align*}
\Delta^{k} (a_{n}b_{n})
=\sum^{k}_{j=0}\dbinom{k}{j}\Delta^{k-j} a_{n+j}\cdot\Delta^{j}b_{n}
\end{align*}
for $ n, k\geq 0.$
This can be shown by induction on $k$ with the
simple identities $\Delta(A_nB_n)=(\Delta A_n)B_n+A_{n+1}\Delta B_n$
and $\binom{k}{j}+\binom{k}{j-1}=\binom{k+1}{j}$ for $1\le j\le k.$
Since $\Delta^{k-j} a_{n+j}\ge0$ and $\Delta^{j}b_{n}\ge 0,$
we obtain $\Delta^k(a_nb_n)\ge0.$
\end{pf}

This result is also claimed by Reza and Zhang \cite[Lemma 1.9]{RZ19}.
According to them, this follows from the fact that $\{a_nb_n\}$ corresponds to
the convolution measure $\mu\diamond\nu$ when
$\{a_n\}$ and $\{b_n\}$ correspond to measures $\mu$ and $\nu,$ respectively.

We are now ready to prove Proposition \ref{proposition1.5}.

\begin{pf}[Proof of Proposition \ref{proposition1.5}]
Let $f_j=h_j+c_j\bar g_j$ for $j=1,2$ and
$$
h_j(z)=\int_0^1\frac{1}{1-tz}d\mu_j(t)
\aand
g_j(z)=\int_0^1\frac{1}{1-tz}d\nu_j(t)
$$
for some Borel probability measures $\mu_j, \nu_j$ for $j=1,2.$
The first assertion immediately follows from the fact that
$(1-s)\mu_1+s\mu_2$ and $(1-s)\nu_1+s\nu_2$ are
Borel probability measures for $0\le s\le 1.$
For the second assertion, we express $f_1\ast f_2$ in the form
$$
(f_1\ast f_2)(z)=(h_1\ast h_2)(z)+c_1c_2\,\overline{(g_1\ast g_2)(z)}.
$$
By Lemma \ref{lem:conv}, we have $h_1\ast h_2, g_1\ast g_2\in\acm.$
Thus the assertion follows.
\end{pf}

\begin{pf}[Proof of Theorem \ref{Theorem 1.6}]
Since $h,g\in \acm$, by Lemma \ref{lem:conv}, $h*g\in\acm$.
Thus, it is easy to see that the function $f$ given in \eqref{equ1.12} belongs to $\HT(c).$

Next we prove the quasiconformality of $f$.
Since $h,g\in \acm$,  $h$ and $g$ can be expressed as
\begin{align}
h(z)=z\sum_{n=0}^{\infty}a_{n}z^{n}, \quad g(z)=z\sum_{n=0}^{\infty}b_{n}z^{n},
\end{align}
for some Hausdorff moment sequences $\{a_{n}\}$ and $\{b_{n}\}.$
Furthermore, there exists a Borel probability measure $\nu$ on $[0,1]$ such that
\begin{align}
b_{n}=\int^{1}_{0}t^{n}d\nu(t), \quad n=0,1,2, \dots.
\end{align}
A simple computation leads to
\begin{align*}
(h\ast g)(z)&=\sum_{n=0}^{\infty}a_{n}b_{n}z^{n+1}
=\sum_{n=0}^{\infty}\left(a_{n}z^{n+1}\int^{1}_{0}t^{n}d\nu(t)\right) \\
&=\int^{1}_{0}\left(\sum_{n=0}^{\infty}a_{n}z^{n+1}t^{n}\right)d\nu(t)
=\int^{1}_{0}\frac{h(tz)}{t}d\nu(t).
\end{align*}
We remark that this property indeed characterizes the generating functions of
Hausdorff moment sequences (see Grinshpan \cite[Theorem 1]{Grin96}).
Thus
\begin{align*} (h*g)'(z)=\int^{1}_{0}h'(tz)d\nu(t).
\end{align*}
and therefore
\begin{align*}
\omega_{f}=\frac{\bar f_z}{f_{z}}
=c\frac{(h*g)'}{h'}=c\int^{1}_{0}\frac{h'(tz)}{h'(z)}d\nu(t).
\end{align*}
By the assumption \eqref{Theorem1.6.1}, we have
\begin{align*}
|\omega_{f}(z)|&= c\left|\int^{1}_{0}\frac{h'(tz)}{h'(z)}d\nu(t)\right| \\
&\leq c  \int^{1}_{0}\left|\frac{h'(tz)}{h'(z)}\right|d\nu(t)\leq cM\int^{1}_{0}d\nu(t)=cM\le k<1.
\end{align*}
In particular, $f$ is locally univalent and thus, by Theorem \ref{Theorem 1.4}, $f$ is univalent on $\D.$
Since $f$ is smooth on $\D$, the inequality
$|\omega_{f}|\le k< 1$ implies that $f$ is $k$-quasiconformal on $\D.$
\end{pf}

\begin{pf}[Proof of Proposition \ref{prop:M}]
First we note that $h'(z)$ vanishes nowhere on $\D$ by assumption.
Let
$$
u(z)=\Re \frac{zh''(z)}{h'(z)}
$$
for $z\in\D.$
Then $u$ is harmonic and $u>-m$ on $\D$ and $u(0)=0.$
If we put $U=(u+m)/m,$ then $U>0$ and $U(0)=1.$
Thus the Harnack inequality implies the inequality $U(z)\ge (1-r)/(1+r)$
for $|z|=r<1.$
Hence,
\begin{equation}\label{eq:u}
u(z)=mU(z)-m\ge \frac{-2mr}{1+r},\quad r=|z|<1.
\end{equation}
Next we set
$$
\psi(s)=\log|h'(sz)|=\Re \log h'(sz),\quad 0\le s\le 1,
$$
for a fixed $z\in\D.$
Then, by \eqref{eq:u},
$$
\psi'(s)=\Re\frac{z h''(sz)}{h'(sz)}=\frac{u(sz)}{s}
\ge\frac{-2ms|z|}{s(1+s|z|)}\ge -2m.
$$
An integration of the above inequality in $t\le s\le 1$ gives us
$$
\log\frac{|h'(z)|}{|h'(tz)|}=\psi(1)-\psi(t)
=\int_t^1\psi'(s)ds\ge -2m(1-t)\ge -2m,
$$
which yields the required inequality.
\end{pf}

\begin{pf}[Proof of Theorem \ref{Theorem 1.8}]
Since the assertion $g/h\in\acm$ and its proof are contained in \cite{RSS09},
we only show the assertion $F:=g'/h'\in\acm.$
Indeed, we will employ the same method as in \cite{RSS09}.

It suffices to check the three conditions in Lemma \ref{lem:char} for $F.$
By the expressions in \eqref{eq:hg}, we have
\begin{align*}
h'(z)=\int^{1}_{0}\frac{\phi(t)}{(1-tz)^{2}}dt
\aand
g'(z)=\int^{1}_{0}\frac{\psi(t)}{(1-tz)^{2}}dt.
\end{align*}
In particular, for a real number $x<1,$ we have $h'(x)\ge0$ and $g'(x)\ge0$
and thus condition (ii) in Lemma \ref{lem:char} is verified.
Condition (i) is clearly satisfied.
The remaining task is to check condition (iii) in the lemma.
Since
$$
\Im F(z)=\frac{\overline{h'(z)}g'(z)-h'(z)\overline{g'(z)}}{2i|h'(z)|^{2}},
$$
we have only to show that $(\overline{h'} g'-h'\overline{g'})/i$ is non-negative on
the upper half-plane $\Im z>0.$
We now compute
\begin{align*} \overline{h'(z)}g'(z)
&=\int^{1}_{0}\frac{\phi(s)}{(1-s\overline{z})^{2}}dt\int^{1}_{0}\frac{\psi(t)}{(1-tz)^{2}}dt\\
&=\int^{1}_{0}\int^{1}_{0}\frac{\phi(s)\psi(t)}{(1-s\overline{z})^{2}(1-tz)^{2}}dsdt\\
&=\int\int_{s\leq t}+\int\int_{t\leq s}\\
&=\int\int_{s\leq t}\left(\frac{\phi(s)\psi(t)}{(1-s\overline{z})^{2}(1-tz)^{2}}
+\frac{\phi(t)\psi(s)}{(1-t\overline{z})^{2}(1-sz)^{2}}\right)dsdt\\
&=\int\int_{s\leq t}\frac{\phi(s)\psi(t)(1-t\overline{z})^{2}(1-sz)^{2}
+\phi(t)\psi(s)(1-s\overline{z})^{2}(1-tz)^{2}}{|1-s\overline{z}|^{4}|1-tz|^{4}}dsdt.
\end{align*}
Taking the complex conjugate, we have similarly
\begin{align*}h'(z) \overline{g'(z)}
=\int\int_{s\leq t}\frac{\phi(s)\psi(t)(1-s\overline{z})^{2}(1-tz)^{2}
+\phi(t)\psi(s)(1-t\overline{z})^{2}(1-sz)^{2}}{|1-s\overline{z}|^{4}|1-tz|^{4}}dsdt.
\end{align*}
We obtain
\begin{align*}
&\overline{h'(z)}g'(z)-h'(z) \overline{g'(z)} \\
=\,&\, 4i\int\int_{s\leq t}\frac{y(t-s)\{1-(s+t)x+str^{2}\}\{\phi(s)\psi(t)-\phi(t)\psi(s)\}}%
{|1-s\overline{z}|^{4}|1-tz|^{4}}dsdt,
\end{align*}
where $z=x+iy, \, r=|z|$.
Since
$$
1-(s+t)x+str^{2}\geq 1-(s+t)x+stx^{2}
=(1-sx)(1-tx)\geq(1-s)(1-t)\geq 0
$$
for $x\in(-\infty,1),$
condition (iii) is now easily confirmed as required.
\end{pf}

\begin{pf}[Proof of Theorem \ref{Theorem 1.9}]
By hypothesis, $\omega_f=\bar f_z/f_z=c g'/h'=cF.$
Note that $h'$ is non-vanishing on $\D$ by Lemma \ref{lem:Wirths}.
The inequality \eqref{eq:F} now leads to $|\omega_f|\le cF(1^-)\le k<1.$
In particular, $f$ is locally univalent and thus Theorem \ref{Theorem 1.4}
implies that $f$ is univalent on $\D.$
We now conclude that $f$ is $k$-quasiconformal on $\D.$
\end{pf}

\section{Examples and applications}

Let us first see an explicit estimate of the constant $M$ in Theorem \ref{Theorem 1.6}.

\begin{example}
Let $h(z)=z/(1-z)$.
Then $h\in \acm$ as we remarked in Introduction.
Since $h'(z)=1/(1-z)^2,$ we have for $z=x+iy$ with fixed $r=|z|<1,$
\begin{align*}
\left|\frac{h'(tz)}{h'(z)}\right|
=\left|\frac{(1-z)^{2}}{(1-tz)^{2}}\right|
=\frac{1-2x+r^2}{1-2tx+t^2r^2}
\le \frac{(1+r)^2}{(1+tr)^2}
<\frac{4}{(1+t)^2}.
\end{align*}
Hence,
$$
\sup_{z\in\D, 0\le t\le 1}\left|\frac{h'(tz)}{h'(z)}\right|=4.
$$
Taking $M=4$ in Theorem \ref{Theorem 1.6},
we know that $f$ given in \eqref{equ1.12} is $4c$-quasiconformal on $\D$
for $0\le c<1/4.$
\end{example}

We next give a simple example to examine the conditions in Theorems
\ref{Theorem 1.4} and \ref{Theorem 1.6}.

\begin{example}\label{ex:qc}
Let $h(z)=z/(1-z)$ as above and $g(z)=z$.
Note that $h,g\in\acm.$
For a positive constant $c<1$ we consider the function
$f(z)=h(z)+c\overline{g(z)}=z/(1-z)+c\bar z$ in $\HT(c).$
Note that $h\ast g=g$ in this case.
The previous example tells us that $f$ is $4c$-quasiconformal on $\D$
for $c<1/4.$
This bound is sharp.
Indeed, the second complex dilatation of $f$ is $\omega_f(z)=cg'(z)/h'(z)=c(1-z)^2$
and thus satisfies $\|\omega_f\|_\infty=4c.$
Moreover, $f$ is not locally univalent on $\D$ for each $c>1/4.$
We will show it.
Let $\gamma$ be the intersection of the circle $|z-1|=1/\sqrt c$ and $\D.$
Note that $\gamma$ is non-empty because $1/\sqrt c<2.$
Points in this arc $\gamma$ may be parametrized as
$z=1+e^{i\theta}/\sqrt c.$
Then
$$
f\big(1+e^{i\theta}/\sqrt c\big)
=\frac{1+e^{i\theta}/\sqrt c}{-e^{i\theta}/\sqrt c}+c(1+e^{-i\theta}/\sqrt c)
=c-1,
$$
which shows that the open arc $\gamma$ shrinks to the one point $c-1.$
Therefore, $f$ is not locally univalent at each point of $\gamma.$
\end{example}

Let us now take a look at polylogarithms.
The polylogarithmic function of order $\alpha$ is defined by
\begin{align*} \Li_{\alpha}(z) =\sum_{n=1}^{\infty}\frac{z^{n}}{n^{\alpha}},
\quad z\in\D,\alpha\geq 0.
\end{align*}
By the well-known representation
$$
\Li_\alpha(z)=\frac{z}{\Gamma(\alpha)}\int_0^1\frac{(-\log t)^{\alpha-1}}{1-tz}dt
$$
for $\alpha>0,$ and $\Li_0(z)=z/(1-z),$
we see that $\Li_\alpha\in\acm$ for $\alpha\ge0.$
Also the relation
\begin{align}\label{polylog 1}
\frac{\Li_{\alpha}}{\Li_{\beta}}\in \cm, \quad 0\leq \alpha\leq\beta,
\end{align}
follows from Theorem \ref{Theorem 1.8} and was already contained
in \cite[Lemma 5.1]{RSS09}.
Lewis \cite{Lewis83} proved that $\Li_\alpha$ maps $\D$ univalently
onto a convex domain for each $\alpha\ge0.$
Furthermore, in \cite{RSS09}, the polylogarithmic function $\Li_\alpha$ is shown to be
universally starlike for $\alpha=0$ and $1\le\alpha$
and universally convex for $\alpha=0,1$ and $2\le \alpha$
and was conjectured to be universally starlike also for $0<\alpha<1$ and
universaly convex for $1<\alpha<2.$
This conjecture was completely proved by
Bakan, Ruscheweyh and Salinas \cite{BRS15}.

We need the following estimate below.
Though this is essentially contained in \cite{RSS09}, we include
a direct proof of it for convenience of the reader.

\begin{lem}\label{lem:T}
Let $F\in\cm$ and $\mu$ be its representing measure on $[0,1].$
Then the following inequalities hold:
$$
\Re F(z)\ge \int_0^1\frac{d\mu(t)}{1+t}\ge \frac12,\quad z\in\D.
$$
\end{lem}

\begin{pf}
By assumption, $F$ is expressed by
$$
F(z)=\int_0^1 \frac{1}{1-tz}d\mu(t), \quad z\in\Lambda.
$$
Letting $z=x+iy$ and $r=|z|<1,$ we compute
$$
\Re F(z)=\int_0^1\frac{1-tx}{1-2tx+t^2r^2}d\mu(t).
$$
Since the function $x\mapsto (1-tx)/(1-2tx+t^2r^2)$ is increasing
in $-r\le x\le r$ for fixed $r$ and $t,$ we have the estimates
$$
\Re F(z)\ge F(-r)
=\int_0^1\frac{1+tr}{1+2tr+t^2r^2}d\mu(t)
=\int_0^1\frac{1}{1+tr}d\mu(t)
\ge \int_0^1\frac{1}{1+t}d\mu(t).
$$
\end{pf}

We apply the above observations to polylogarithms to have the following.

\begin{thm}\label{Theorem 1.7}
Let $\alpha, \beta\ge 1$ and $c$ be a non-negative real constant
and set $f=\Li_\alpha+c\overline{\Li_\beta}.$
\begin{enumerate}
\item[(i)]
$f$ is $k$-quasiconformal on $\D$ when $\alpha\le\beta$ and $2c\le k<1;$
\item[(ii)]
$f$ is $k$-quasiconformal on $\D$ when $2<\beta\le\alpha$ and
$c\,\zeta(\beta-1)/\zeta(\alpha-1)\le k<1.$
\end{enumerate}
\end{thm}

Here, $\zeta(s)$ denotes the Riemann zeta function $s\mapsto
\displaystyle\sum_{n=1}^\infty n^{-s}.$
Recall that $\zeta(s)<+\infty$ for $s\in(1,+\infty).$

\begin{pf}
Put $h=\Li_\alpha$ and $g=\Li_\beta$ for brevity.
Note that they are univalent on $\D$ by Lewis' theorem \cite{Lewis83}.
Since $zh'(z)=\Li_{\alpha-1}(z),$ $zg'(z)=\Li_{\beta-1}(z)$ we have
$$
F:=\frac{g'}{h'}=\frac{zg'}{zh'}=\frac{\Li_{\beta-1}}{\Li_{\alpha-1}}.
$$
First assume that $\alpha\le\beta.$
Then $G:=1/F\in\cm$ by \eqref{polylog 1}.
Hence, Lemma \ref{lem:T} implies $|G(z)|\ge \Re G(z)\ge 1/2$
for $z\in\D.$
In view of the form of $f,$ we now estimate
$$
|\omega_f|
=\left|\frac{\bar f_z}{f_z}\right|
=\left|\frac{cg'}{h'}\right|
=\frac{c}{|G|}\le 2c\le k<1
$$
on $\D.$
In particular, $f$ is locally univalent and thus, by Theorem \ref{Theorem 1.4}, $f$ is univalent on $\D.$
It is now clear that $f$ is $k$-quasiconformal on $\D.$

Next assume that $\alpha\ge\beta.$
Then $F\in\cm$ by \eqref{polylog 1}.
If $\beta>2,$ we have $F(1^-)=\Li_{\beta-1}(1^-)/\Li_{\alpha-1}(1^-)
=\zeta(\beta-1)/\zeta(\alpha-1)<+\infty.$
Now the assertion follows from Theorem \ref{Theorem 1.9}.
\end{pf}

We remark that we have $F(1^-)=+\infty$ when $0\le\beta\le 2$ and $\beta<\alpha.$

\medskip

Finally, we apply our results to hypergeometric functions.
We recall the definition of the hypergeometric function $\Gauss(a,b;c;z)$:
$$
\Gauss(a,b;c;z)=1+\sum_{n=1}^\infty \frac{(a)_n(b)_n}{(c)_n n!}z^n,\quad |z|<1,
$$
where $(a)_n$ is the Pochhammer symbol; namely, $(a)_n=a(a+1)\cdots(a+n-1)$ for
$n\ge 1$ and $(a)_0=1.$
Here, $a,b$ and $c$ are (possibly complex) parameters with $c\ne 0,-1,-2,\dots.$
Geometric properties such as starlikeness and convexity of $\Gauss(a,b;c;z)$ and the
shifted one $z\Gauss(a,b;c;z)$ were studied by many authors
(see \cite{Kus07}, \cite{SW17geo}, \cite{SW17sp} and references therein).
In particular,  in connection with the class $\cm,$ some observations on the
hypergeometric functions were made in \cite{RSS09}.
The formula
\begin{equation}\label{eq:Gauss}
\Gauss(a,b;c;1^-)=\frac{\Gamma(c)\Gamma(c-a-b)}{\Gamma(c-a)\Gamma(c-b)},
\quad \Re(c-a-b)>0,
\end{equation}
is due to Gauss.
It should also be note that 
the derivatve formula $\frac{d}{dz}\Gauss(a,b;c;z)=\frac{ab}c\Gauss(a+1,b+1;c+1;z)$ holds.
The well-known Euler representation formula
$$
\Gauss(a,b;c;z)=\frac{\Gamma(c)}{\Gamma(a)\Gamma(c-a)}\int_0^1(1-tz)^{-b}t^{a-1}(1-t)^{c-a-1}dt
$$
for $\Re c>\Re a>0$ implies that the function
$$
L_{a,c}(z)=z\Gauss(a,1;c;z)
$$
belongs to the class $\acm$
for real parameters $c>a>0$ with the representing measure $\mu$ given by
$$
d\mu(t)=\frac{\Gamma(c)}{\Gamma(a)\Gamma(c-a)}t^{a-1}(1-t)^{c-a-1}dt.
$$
We note that $L_{a,c}(z)$ is univalent on the half-plane $\Re z<1$
by Lemma \ref{lem:Wirths} and therefore $L_{a,c}'$ is non-vanishing there.
The convolution $f\ast L_{a,c}$ with analytic functions $f$
is often called the Carlson-Shaffer operator and studied by
many authors.

As a simple application of Theorem \ref{Theorem 1.8}, we have the following result.

\begin{lem}\label{lem:L}
Let $a,c, a', c'$ be real constants with $c>a>0$ and $c'>a'>0.$
If $a'\ge a$ and if $c-a\ge c'-a',$ then the functions $L_{a',c'}/L_{a,c}$
and $L_{a',c'}'/L_{a,c}'$ both belong to the class $\cm.$
\end{lem}

In a similar way to Theorem \ref{Theorem 1.7}, we finally obtain the following.

\begin{thm}
Let $a,c, a', c'$ be real constants with $c>a>0$ and $c'>a'>0.$
For a non-negative real constant $b,$ set $f=L_{a,c}+b\,\overline{L_{a',c'}}.$
\begin{enumerate}
\item[(i)]
$f$ is $k$-quasiconformal on $\D$ when $a\ge a',~ c-a\le c'-a'$
and $2b\le k<1;$
\item[(ii)]
$f$ is $k$-quasiconformal on $\D$ when $a'\ge a,~ 2<c'-a'\le c-a$ and
$bM\,\le k<1,$ where
$$
M=\frac{(c'-1)(c'-2)(c-a-1)(c-a-2)}{(c-1)(c-2)(c'-a'-1)(c'-a'-2)}.
$$
\end{enumerate}
\end{thm}

\begin{pf}
Put $h=L_{a,c}$ and $g=L_{a',c'}$ for brevity.
First assume that $a\ge a',~ c-a\le c'-a'.$
Then by the previous lemma, $G=h'/g'\in\cm.$
Therefore, Lemma \ref{lem:T} implies $|G|\ge \Re G\ge 1/2$ on $\D.$
We estimate as before
$$
|\omega_f|
=\left|\frac{bg'}{h'}\right|
=\frac{b}{|G|}\le 2b\le k<1
$$
on $\D$ and thus conclude that $f$ is $k$-quasiconformal on $\D.$

Next assume that $a'\ge a$ and $2<c'-a'\le c-a.$
Then, by Lemma \ref{lem:L}, we have $F=g'/h'\in\cm.$
Note that $h'(z)=\Gauss(a,1;c;z)+(a/c)z\Gauss(a+1,2;c+1;z).$
By \eqref{eq:Gauss} and the basic identity $\Gamma(x+1)=x\Gamma(x),$
\begin{align*}
h'(1^-)&=\frac{\Gamma(c)\Gamma(c-a-1)}{\Gamma(c-a)\Gamma(c-1)}
+\frac ac\cdot \frac{\Gamma(c+1)\Gamma(c-a-2)}{\Gamma(c-a)\Gamma(c-1)} \\
&=\frac{c-1}{c-a-1}+\frac ac\cdot \frac{c(c-1)}{(c-a-1)(c-a-2)}
=\frac{(c-1)(c-2)}{(c-a-1)(c-a-2)}.
\end{align*}
Similarly, we have
$$
g'(1^-)=\frac{(c'-1)(c'-2)}{(c'-a'-1)(c'-a'-2)}.
$$
Hence, $F(1^-)=M<+\infty.$
Now the assertion follows from Theorem \ref{Theorem 1.9}.
\end{pf}

\def\cprime{$'$} \def\cprime{$'$} \def\cprime{$'$}
\providecommand{\bysame}{\leavevmode\hbox to3em{\hrulefill}\thinspace}
\providecommand{\MR}{\relax\ifhmode\unskip\space\fi MR }
\providecommand{\MRhref}[2]{%
  \href{http://www.ams.org/mathscinet-getitem?mr=#1}{#2}
}
\providecommand{\href}[2]{#2}

\end{document}